\documentclass[10pt]{article}   
\usepackage{amsthm, amsmath}
\usepackage{amssymb}   
\theoremstyle{plain}

\newtheorem{theorem}{Theorem}[section]
\newtheorem{lemma}[theorem]{Lemma}
\newtheorem{proposition}[theorem]{Proposition}
\newtheorem{corollary}[theorem]{Corollary}

\theoremstyle{definition}

\newtheorem{definition}[theorem]{Definition}
\newtheorem{example}[theorem]{Example}

\theoremstyle{remark}

\newtheorem{remark}[theorem]{Remark}

\title{COMPLEXITY OF VILLAMAYOR'S ALGORITHM \\ IN THE NON EXCEPTIONAL MONOMIAL CASE}
\author{Rocio Blanco\footnote{Research partially supported by an F.P.U. Fellowship, Spanish Ministry of Education and Culture, AP2002-0009.}}
\date{}

\begin{document}
\maketitle

\begin{abstract}
We study monomial ideals, always locally given by a monomial, like
a reasonable first step to estimate in general the number of monoidal transformations of Villamayor's
algorithm of resolution of singularities. The resolution of a monomial
ideal $<X_1^{a_1}\cdot \ldots \cdot X_n^{a_n}>$ is interesting due
to its identification with the particular toric problem $<Z^c-
X_1^{a_1}\cdot \ldots \cdot X_n^{a_n}>$.

In the special case, when all the exponents $a_i$ are greater than or
equal to the critical value $c$, we construct the largest branch of the resolution tree which provides 
an upper bound involving partial sums of Catalan numbers. This case will be
called ``minimal codimensional case''. Partial sums of Catalan numbers (starting
$1,2,5,\ldots$) are $1,3,8,22,\ldots$ These partial sums are well known in Combinatorics and count the number of paths starting from the root in all ordered trees with $n+1$ edges. Catalan numbers appear in many
combinatorial problems, counting the number of ways to insert $n$
pairs of parenthesis in a word of $n+1$ letters, plane trees with
$n+1$ vertices, $\ldots $, etc.

The non minimal case, when there exists some exponent $a_{i_0}$ smaller than $c$, will be called ``case of higher codimension''. In this case, still unresolved, we give an example to state the foremost troubles.

Computation of examples has been helpful in both cases to study the behaviour of the resolution invariant. Computations have been made in Singular (see \cite{sing}) using the \emph{desing} package by G. Bodn\'ar and J. Schicho, see \cite{lib}. 
\end{abstract}



\section{Introduction}

The existence of resolution of singularities in arbitrary dimension over a field of characteristic zero was solved by Hironaka in his famous paper \cite{existence}. Later on, different constructive proofs have been given, among others, by Villamayor \cite{villa}, Bierstone-Milman \cite{bm}, Encinas-Villamayor \cite{cour}, Encinas-Hauser \cite{strong} and Wodarczyk \cite{Wo}. 

This paper is devoted to study the complexity of Villamayor's algorithm of resolution of singularities. This algorithm appears originally in \cite{villa} and we will use the presentation given in \cite{cour}. In this paper, the authors introduce a class of objects called \emph{basic objects} $B=(W,(J,c),E)$ where $W$ is a regular ambient space over a field $k$ of characteristic zero, $J\subset \mathcal{O}_W$ is a sheaf of ideals, $c$ is an integer and $E$ is a set of smooth hypersurfaces in $W$ having only normal crossings. That is, they consider the ideal $J$ together with a positive integer $c$, or \emph{critical value} defining the \emph{singular locus} $Sing(J,c)=\{\xi \in W|\ ord_{\xi}(J)\geq c \}$, where $ord_{\xi}(J)$ is the \emph{order} of $J$ in a point $\xi$.

Let $W\stackrel{\pi}{\leftarrow} W'$ be the monoidal transformation with center $\mathcal{Z}\subset Sing(J,c)$, $\pi^{-1}(\mathcal{Z})=Y'$ is the exceptional divisor. Let $\xi$ be the generic point of $\mathcal{Z}$, $ord_{\xi}(J)=\theta$, the total transform of $J$ in $W'$ satisfies $J\mathcal{O}_{W'}=I(Y')^{\theta}\cdot J^{\curlyvee}$ where $J^{\curlyvee}$ is the \emph{weak} transform of $J$, (see \cite{cour} for details). 

A \emph{transformation} of a basic object $(W,(J,c),E)\leftarrow (W',(J',c),E')$ is defined by a monoidal transformation $W\stackrel{\pi}{\leftarrow} W'$ and defining $J'=I(Y')^{\theta-c}\cdot J^{\curlyvee}$, the \emph{controlled transform} of $J$. 

A sequence of transformations of basic objects \begin{equation} (W,(J,c),E)\leftarrow (W^{(1)},(J^{(1)},c),E^{(1)})\leftarrow \cdots \leftarrow (W^{(N)},(J^{(N)},c),E^{(N)}) \label{resolution} \end{equation} is a \emph{resolution} of $(W,(J,c),E)$ if $Sing(J^{(N)},c)=\emptyset$. 

\begin{remark} Superscripts $^{(k)}$ in basic objects will denote the $k$-stage of the resolution process. Subscripts $_i$ will always denote the dimension of the ambient space $W_i^{(k)}$. 
\end{remark}

Villamayor's algorithm provides a \emph{log-resolution} in characteristic zero. A log-resolution of $J$ is a sequence of monoidal transformations at regular centers as (\ref{resolution}) such that each center has normal crossings with the exceptional divisors $E^{(i)}$, and the total transform of $J$ in $W^{(N)}$ is of the form  $$J\mathcal{O}_{W^{(N)}}=I(H_1)^{b_1}\cdot \ldots \cdot I(H_N)^{b_N}$$ with $b_i\in \mathbb{N}$ for all $1\leq i \leq N$  and $E^{(N)}=\{H_1,\ldots,H_N\}$. 

In \cite{cour} it is shown that algorithmic principalization of ideals reduces to algorithmic resolution of basic objects. That is, starting with $c\!=\!max\ ord(J)$, the maximal order of $J$, we obtain a resolution of $(W,(J,c),E)$ as (\ref{resolution}). At this step $max\ ord(J^{(N)})\!=\!c^{(N)}\!<c$. If $c^{(N)}\!>1$, we continue resolving $(W^{(N)},(J^{(N)},c^{(N)}),E^{(N)})$ and so on, until have $max\ ord(J^{(\mathcal{N})})\!=\!c^{(\mathcal{N})}\!=1$. Finally, a resolution of $(W^{(\mathcal{N})},(J^{(\mathcal{N})},1),E^{(\mathcal{N})})$ provides a log-resolution of $J^{(\mathcal{N})}$, and therefore a log-resolution of $J$. 

In \cite{cour} it is also shown that algorithmic principalization of ideals leads to embe\-dded desingulari\-zation of varieties. That is, given a closed subscheme $X\subset W$, the algorithmic principalization of the ideal $I(X)$ provides an embedded desingularization of $X$. See also \cite{EncinasVillamayor2003} for more details.

A key point in the definition of the algorithm is to use induction on the dimension of the ambient space $W$ to define an upper-semi-continuous function $t$. The set of points where this function attains its maximal value, $\underline{Max}\ t$, is a regular closed set, and defines a regular center for the next monoidal transformation. 

A resolution of the basic object $(W,(J,c),E)$ is achieved by a sequence of monoidal transformations as in (\ref{resolution}), with centers $\underline{Max}\ t^{(k)}$ for $0\leq k\leq N-1$. That is, the sequence of 
monoidal transformations is defined by taking successively the center defined by the upper-semi-continuous function. The algorithm stops at some stage because the maximal value of the function $t$ drops after monoidal transformations, that is, $max\ t^{(0)}>max\ t^{(1)}>\ldots >max\ t^{(N-1)}.$ 

This function $t$ will be the resolution invariant. We shall work with the invariant defined in \cite{cour},
using the language of mobiles developed in \cite{strong}. We
remind briefly the main notions. 

Let $J\subset \mathcal{O}_W$ be an ideal defining a singular algebraic set $X\subset W$. The ideal $J$ factors into $J=M\cdot I$, with $M$ the ideal defining a normal crossing divisors, and $I$ some ideal
still unresolved. 

By induction on the dimension of $W$, we will
have this decomposition at every dimension from $n$ to $1$, that
is $J_i=M_i\cdot I_i$, for $n\geq i \geq 1$, are defined in local flags $W_n\supseteq W_{n-1}\supseteq \cdots \supseteq W_i\supseteq  \cdots \supseteq W_1$, where each $J_i,M_i,I_i \in \mathcal{O}_{W_i}$ are in dimension $i$. There is a
critical value $c_{i+1}$ at each dimension $i$, ($c_{n+1}=c$), see \cite{strong} for details. All the basic objects $(W_i,(J_i,c_{i+1}),E_i)$, for $n\geq i \geq 1$, will be resolved during the process of the algorithm. 

Let $E$ be the exceptional divisor of previous monoidal transformations, and consider $E=\cup_{i=1}^nE_i$ where $E_i$ applies to dimension $i$. Obviously, we start with $E=\emptyset$.

For any point $\xi \in Sing(J,c)$, the function $t$ will
have $n$ coordinates, with lexicographical order, and it will be
one of the following three types: \begin{equation} \label{invt} \hspace*{-0.2cm} \begin{array}{ll} (a) &
t(\xi)=(t_n(\xi),t_{n-1}(\xi),\ldots,t_{n-r}(\xi),\ \infty,\
\infty,\ldots,\infty)
\\ (b) & t(\xi)=(t_n(\xi),t_{n-1}(\xi),\ldots,t_{n-r}(\xi),
\Gamma(\xi),\infty,\ldots,\infty)
\\ (c) & t(\xi)=(t_n(\xi),t_{n-1}(\xi),\ldots,t_{n-r}(\xi),\ldots
\ldots \ldots,t_1(\xi)) \end{array}\  \text{ with }
t_i=\left[\frac{\theta_i}{c_{i+1}},m_i\right] \end{equation} where
$\theta_i=ord_{\xi}(I_i)\ ,$ $m_i$ is the number of exceptional
divisors in $E_i$, and $\Gamma$ is the resolution function corresponding to the so-called \emph{monomial case}, following the notation of \cite{cour}, pages $165-166$. We will recall the definition of $\Gamma$ in equation (\ref{gamma}).

For simplicity, let assume that we start with a polynomial ring, $W=Spec(k[X_1,\ldots,X_n])$. In $\mathcal{O}_W=k[X_1,\ldots,X_n]$ the ideal $J$ is locally given by a monomial with respect to a regular system of parameters $$J=<X_1^{a_1}\cdot \ldots \cdot X_n^{a_n}>\subset \mathcal{O}_W \text{ with } a_i\in \mathbb{N}, \text{ for } i=1,\ldots,n.$$

Note that, in this situation, the center of the next monoidal transformation is combinatorial, it is a linear combination of $X_1,\ldots,X_n$. And this is also true after monoidal transformations, since Villamayor's algorithm applied to a monomial ideal provides always combinatorial centers, and after a monoidal transformation in a combinatorial center we obtain again a monomial ideal.    

So, at any stage of the resolution process, $W=\cup_i U_i$, where $U_i\cong\mathbb{A}^n_k$. Thereafter, we shall work locally, so we will assume that $W$ is an affine space.

To resolve the toric hypersurface $\{f=0\}=\{Z^c-X_1^{a_1}\cdot \ldots \cdot X_n^{a_n}=0\}$ we note that its singular locus $Sing(<f>,c)$ is always included in $\{Z=0\}$, so we argue by induction on the dimension and reduce to the case
where the corresponding ideal $J$ is of the form
\begin{equation} \label{jota} J=<X_1^{a_1}\cdot \ldots \cdot
X_n^{a_n}>\subset \mathcal{O}_W \text{ with }\ 1\leq a_1\leq a_2\leq \ldots \leq a_n,\
\sum_{i=1}^n a_i=d,\ d\geq c,\end{equation} where $c$ is
the critical value. If $a_i=0$ for some $i$, then we may assume $dim(W)<n$. 

After a monoidal transformation, we always consider the controlled transform of $J$ with respect to $c$, $J'=I(Y')^{-c}\cdot J^{*}$ where $J^{*}$ is the total transform of $J$ and $Y'$ denotes the new exceptional divisor. For the toric problem $J=<Z^c-X_1^{a_1}\cdot \ldots \cdot X_n^{a_n}>$, taking the origin as center of the next monoidal transformation, at the $i$-th chart: $$J^*=<Z^c\cdot X_i^{c} - X_1^{a_1}\cdots X_i^{d} \cdots X_n^{a_n}>=<X_i^{c}\cdot(Z^c-X_1^{a_1}\cdots X_i^{d-c} \cdots X_n^{a_n})>,$$ and we can only factorize $c$ times the exceptional divisor. 

\begin{remark} We will denote as {\it $i$-th chart} the chart where we divide by $X_i$. When the center of the monoidal transformation is the origin, this monoidal transformation is expressed: 
$$\begin{array}{ccc} k[Z,X_1,\ldots,X_n] & \rightarrow & k[Z,X_1,\ldots,X_n,\frac{Z}{X_i},\frac{X_1}{X_i},\ldots,\frac{X_{i-1}}{X_i},\frac{X_{i+1}}{X_i},\ldots,\frac{X_n}{X_i}] \\ Z & \rightarrow & \frac{Z}{X_i} \\ X_i & \rightarrow & X_i \\ X_j & \rightarrow & \hspace*{1.5cm} \frac{X_j}{X_i} \text{ for } j\neq i \\
\end{array} $$  where $k[Z,X_1,\ldots,X_n,\frac{Z}{X_i},\frac{X_1}{X_i},\ldots,\frac{X_{i-1}}{X_i},\frac{X_{i+1}}{X_i},\ldots,\frac{X_n}{X_i}]\cong \\ k[\frac{Z}{X_i},\frac{X_1}{X_i},\ldots,\frac{X_{i-1}}{X_i},X_i,\frac{X_{i+1}}{X_i},\ldots,\frac{X_n}{X_i}]$. For simplicity, we will denote each $\frac{X_j}{X_i}$ again as $X_j$, and $\frac{Z}{X_i}$ as $Z$.
\end{remark} 

So we will apply the resolution algorithm to the basic object $(W,(J,c),\emptyset)$ for $J=<X_1^{a_1}\cdot \ldots \cdot X_n^{a_n}>$, which is already a monomial ideal, but it is not necessarily supported on the exceptional divisors. 

\section{Monomial case (exceptional monomial)}

The \emph{monomial case} is a special case in which J is a ``monomial ideal" given locally by a monomial that can be expressed in terms of the exceptional divisors. This case arises after several monoidal transformations.

This means that we have a basic object $(W,(J,c),E)$ where $J$ is locally defined by one monomial supported on the hypersurfaces in $E$. In this case, the ideal $J$ factors into $J=M\cdot I$ with $J=M$ and $I=1$. We can also call it {\bf exceptional monomial}.

\begin{theorem} \label{mon} Let $J\subset \mathcal{O}_W$ be a monomial ideal as in equation (\ref{jota}). Let
$E=\{H_1,\ldots,H_n\}$, with $H_i=V(X_i)$, be a normal crossing divisor. \\ Then an upper bound for
the number of monoidal transformations to resolve $(W,(J,c),E)$ is given by $$\frac{d - c +
gcd (a_1,\ldots,a_n,c)}{gcd (a_1,\ldots,a_n,c)}.$$
\end{theorem}

\begin{proof} We may assume that the greatest common divisor of the
exponents $a_i$ and the critical value $c$ is equal to $1$, because
both the simplified problem and the original problem have the
same singular locus. That is, if $gcd (a_1,\ldots,a_n,c)=k$ then $d=k\cdot d_1$, $c=k\cdot c_1$, $a_i=k\cdot b_i$ for all $1\leq i\leq n$ and $gcd (b_1,\ldots,b_n,c_1)=1$. The ideal $J$ can be written as $J=(J_1)^k$ where $J_1=<X_1^{b_1}\cdot \ldots \cdot X_n^{b_n}>$ therefore $$Sing(J,c)=\{\xi \in X|\ ord_{\xi}((J_1)^k)\geq k\cdot c_1 \}=\{\xi \in X|\ ord_{\xi}(J_1)\geq c_1 \}=Sing(J_1,c_1),$$ where $X$ is the algebraic set defined by $J$.

For a point $\xi \in \mathbb{A}^n_k$, $\Gamma(\xi)=(-\Gamma_1(\xi),\Gamma_2(\xi),\Gamma_3(\xi))$ where \begin{equation} \label{gamma} \vspace*{0.15cm} \begin{array}{l} \Gamma_1(\xi)=\min\{p\ |\ \exists\ i_1,\ldots,i_p, a_{i_1}(\xi)+\cdots+a_{i_p}(\xi)\geq c,\ \xi\in H_{i_1}\cap\cdots \cap H_{i_p}\}, \vspace*{0.15cm} \\ \Gamma_2(\xi)=\max\left\{\frac{a_{i_1}(\xi)+\cdots+a_{i_p}(\xi)}{c}\ |\ {\scriptstyle p=\Gamma_1(\xi),\ a_{i_1}(\xi)+\cdots+a_{i_p}(\xi)\geq c,\ \xi\in H_{i_1}\cap\cdots \cap H_{i_p}}\right\},\vspace*{0.15cm} \\  \Gamma_3(\xi)=\max\{(i_1,\ldots,i_p,0,\ldots,0)\in \mathbb{Z}^n\ |\ {\scriptstyle \Gamma_2(\xi)=\frac{a_{i_1}(\xi)+\cdots+a_{i_p}(\xi)}{c},\ \xi\in H_{i_1}\cap\cdots \cap H_{i_p}}\} \vspace*{0.1cm} \end{array} \end{equation} with lexicographical order in $\mathbb{Z}^n$.

The center $\mathcal{Z}$ of the next monoidal transformation is given by the set of points where $\Gamma$ attains its maximal value. It is easy to see that $\mathcal{Z}=\cap_{i=n-(r-1)}^n H_i$. 

So at the $j$-th chart, the exponent of $X_j$ after the monoidal transformation is $(\sum_{i=n-r+1}^n a_i)-c$ and   $$\left(\sum_{i=n-r+1}^n a_i\right)-c< \min_{n-r+1\leq i \leq n}a_i=a_{n-r+1}$$ because $\sum_{i=n-r+2}^{n} a_i<c$ by construction of the center $\mathcal{Z}$.

This shows that the order of the ideal drops after each monoidal transformation by at least 
one, so in the worst case, we need $d-(c-1)$ monoidal transformations to obtain an
order lower than $c$.
\end{proof}

\begin{remark} Note that it is necessary to consider the monomial case. On one hand, this case may appear in dimension $n$, and also in lower dimensions, $n-1,\ldots,1$, when we resolve any basic object $(W,(J,c),E)$ (where $J$ is any ideal). So we need to resolve the monomial case in order to obtain a resolution of the original basic object $(W,(J,c),E)$. 

On the other hand, the algorithm of resolution leads to the monomial case, since given any ideal $J$, the algorithm provides a log-resolution of $J$. And it is necessary to continue to a resolution within the monomial case. 
\end{remark}

\begin{remark}
The bound in theorem \ref{mon} is reached only for the following values of $c$: $$1,\
a_n+ \ldots + a_j+1 \text{ for } n\geq j\geq 2,\ d.$$ For these values of $c$, the order of the ideal drops after each monoidal transformation exactly by one: 
\begin{itemize}
	\item If $c=1$, the monoidal transformation is an isomorphism. The exponent of $X_n$ after the monoidal transformation is $a_n-1$.
	\item If $c=a_n+ \ldots + a_j+1$, for $n\geq j\geq 2$, the center of the monoidal transformation is $\mathcal{Z}=\cap_{i=j-1}^n H_i$. At the $l$-th chart, for $n\geq l\geq j-1$, the exponent of $X_l$ after the monoidal transformation is $(\sum_{i=j-1}^n a_i)-c=(\sum_{i=j-1}^n a_i)-(\sum_{i=j}^n a_i)-1=a_{j-1}-1$.
	
In particular, at the $(j-1)$-th chart, the exponent of $X_{j-1}$ after the monoidal transformation has droped exactly by one.  
	\item If $c=d$, we finish after only one monoidal transformation. 
\end{itemize} 
\end{remark}

\begin{remark} If $gcd (a_1,\ldots,a_n,c)=k>1$, then the bound of the theorem \ref{mon}
is $(d-c+k)/ k$. \\ As $(d-c+k)/ k<d-c+1,$ we can use in practice the bound for the case $gcd
(a_1,\ldots,a_n,c)=1$.
\end{remark}

\section{Case of one monomial}

To construct an upper bound for the number of monoidal transformations needed to resolve the basic object
$(W,(J,c),E=\emptyset)$, where $J$ is locally defined by a unique monomial, we estimate the number of
monoidal transformations needed to obtain $(W',(J',c),E')$, a transformation of the original basic object, with $J'=M'$ (an exceptional monomial), and then apply theorem \ref{mon}. In order to use theorem \ref{mon}, we need an estimation of the
order of $M'$. This estimation will be valid at any stage of the resolution process.

\begin{lemma} \label{grad} Let $(W,(J,c),\emptyset)$ be a basic object where $J$ is a monomial ideal as in equation (\ref{jota}). Let $J=M\cdot I$ be the factorization of $J$, where $M=1$, because of $E=\emptyset$, and $J=I$. After $N$ monoidal transformations we have $(W^{(N)},(J^{(N)},c),E^{(N)})$. Let $\xi\in W^{(N)}$ be a point. Then $$\ ord_{\xi}(M^{(N)}) \leq (2^N-1)(d-c)$$
where $ord_{\xi}(M^{(N)})$ denotes the order at $\xi$ of $M^{(N)}$, the (exceptional)
monomial part of $J^{(N)}$.
\end{lemma}

\begin{proof} It follows by induction on $N$: 
\begin{itemize}
\item that if $N=1, \  ord_{\xi}(M^{(1)})=d-c.$ \\ At the beginning, the first center defined by this algorithm is always the origin, so at the $i$-th chart: $$J^{(1)}=M^{(1)} \cdot I^{(1)}=<X_i^{d-c}>\cdot <X_1^{a_1}\cdot \stackrel{\widehat{i}}{\ldots} \cdot X_n^{a_n}>$$ with $E^{(1)}=\{H_i\}$, where $H_i=V(X_i)$.   
\item We assume that the result holds for $N=m-1$. $$J^{(m-1)}= M^{(m-1)} \cdot I^{(m-1)}= <X_{i_1}^{b_1}\cdots
X_{i_s}^{b_s}>\cdot <X_{i_{s+1}}^{a_{i_{s+1}}}\cdots X_{i_n}^{a_{i_n}}>$$ with $\sum_{i=1}^s b_i=d'.$ By inductive hypothesis, after $m-1$ monoidal transformations, the maximal order $d'$ of the (exceptional) monomial part $M^{(m-1)}$ satisfies $$d'\leq (2^{m-1}-1)(d-c).$$
For $N=m$, there are two possibilities:
\begin{enumerate}
\item If $max\ ord(I^{(m-1)})=\sum_{j=s+1}^n a_{i_j}\geq c$ then the center of the next monoidal transformation contains only variables appearing in $I^{(m-1)}$. 
\item If $max\ ord(I^{(m-1)})=\sum_{j=s+1}^n a_{i_j}<c$ then the center of the next monoidal transformation contains variables appearing in $I^{(m-1)}$ and also variables appearing in $M^{(m-1)}$. 
\end{enumerate}
\begin{itemize}
\item[\underline{Case $1$}:] \ \\ If the center of the monoidal transformation is as small as possible, that is $\mathcal{Z}\!=\!\cap_{j=s+1}^n V(X_{i_j})$, at the $i_l$-th chart,
$$J^{(m)}\!=\!M^{(m)}\cdot I^{(m)}\!=<\!X_{i_1}^{b_1}\!\cdots\!X_{i_s}^{b_s}\cdot
X_{i_l}^{d-\sum_{j=1}^s a_{i_j}-c}\!>\!\cdot\!
<\!X_{i_{s+1}}^{a_{i_{s+1}}}\!\stackrel{\widehat{i_l}}{\cdots}\!X_{i_n}^{a_{i_n}}\!>\!.$$ The exponent of $X_{i_l}$, $\ d-\sum_{j=1}^s a_{i_j}-c=\sum_{j=s+1}^n a_{i_j}-c \ $ is as big as possible, so this is the worst case, because the increase in the order of the exceptional monomial part after the monoidal transformation will be greater than that for another centers. 

The highest order of $M^{(m)}$ is $$\sum_{i=1}^s b_i + d-
\sum_{j=1}^s a_{i_j}-c =d'+d-c- \sum_{j=1}^s a_{i_j} \leq
d'+d-c,$$ so by inductive hypothesis $$d'+d-c \leq (2^{m-1}-1)(d-c)+d-c=2^{m-1}(d-c)\leq
(2^m-1)(d-c).$$
\item[\underline{Case $2$}:]
\begin{enumerate}
\item[-] At the $i_j$-th chart, for $1\leq j\leq s$ 
$$ J^{(m)}=M^{(m)}\cdot I^{(m)}=<X_{i_1}^{b_1}
\stackrel{\widehat{i_j}}{\cdots} X_{i_s}^{b_s}\cdot
X_{i_j}^{\square}>\cdot <X_{i_{s+1}}^{a_{i_{s+1}}}
\cdots X_{i_n}^{a_{i_n}}>.$$
\item[-] At the $i_l$-th chart, for $s+1\leq l\leq n$
$$J^{(m)}=M^{(m)}\cdot I^{(m)}=
<X_{i_1}^{b_1}\cdots X_{i_s}^{b_s}\cdot
X_{i_l}^{\vartriangle}>\cdot <X_{i_{s+1}}^{a_{i_{s+1}}}
\stackrel{\widehat{i_l}}{\cdots} X_{i_n}^{a_{i_n}}>.$$
\end{enumerate}
As above, if we are in the worst case, when the center of the monoidal transformation is as small as possible, that is, the center is a point,
$$\square=\vartriangle= d'+d-\sum_{j=1}^s a_{i_j}-c\ .$$ Therefore in both cases the highest order of $M^{(m)}$ satisfies $$\leq 2d'+d-c\ \leq \ 2(2^{m-1}-1)(d-c)+d-c=(2^m-1)(d-c)\ .$$
\end{itemize}
\end{itemize} 
\end{proof}

\begin{remark} Due to its general character, this bound is
large and far from being optimal.
\end{remark}

\begin{remark} \label{constructM} The ideals $M_i$ are supported on a normal crossing divisors $D_i$. Recall that their transformations after monoidal transformations, in the neighbourhood of a point $\xi\in W_i$, are $$D_i'=\left\{ \begin{array}{ll} D_i^*+(\theta_i-c_{i+1})\cdot Y' &
\text{ \scriptsize if } {\scriptstyle
(t_n'(\xi'),\ldots,t_{i+1}'(\xi')=(t_n(\xi),\ldots,t_{i+1}(\xi))}  \\ \emptyset & \text{ \scriptsize in other case} \end{array} \right. ,\ n\geq i\geq 1,$$ 
\hspace*{2.5 cm} (${\scriptstyle D_n'=D_n^*+(\theta_n-c)\cdot Y'}$ {\scriptsize always}) \vspace*{2 mm} \\ where $D_i^*$ denotes the pull-back of $D_i$ by the monoidal transformation $\pi$, $Y'$ denotes the new exceptional divisor, the point $\xi'\in W_i'$ satisfies $\pi(\xi')=\xi$, $\theta_i=ord_{\xi}(I_i)$ and $c_{i+1}$ is the corresponding critical value. 
\end{remark}

In what follows we will define the ideals $J_{i-1}$, $n\geq i> 1$. We need some auxiliary definitions: the companion ideals $P_i$ and the composition ideals $K_i$, see \cite{strong} for details. \\ We construct the companion ideals to ensure that $Sing(P_i,\theta_i)\subset Sing(J_i,c_{i+1})$,
\begin{equation} \label{pes}
P_i= \left\{\begin{array}{ll} I_i & \text{ if }\ \theta_i\geq
c_{i+1} \\ I_i+M_i^{\frac{\theta_i}{c_{i+1}-\theta_i}} & \text{
if }\ 0< \theta_i< c_{i+1} \end{array}\right.
\end{equation} where $\xi\in \mathbb{A}^n_k$ is a point, $\theta_i=ord_{\xi}(I_i)$ and $c_{i+1}$ is the corresponding critical value. 

Let $J_i=M_i\cdot I_i$ be the factorization of an ideal $J_i$ in $W_i$, where $M_i,I_i$ are ideals in $W_i$ in the neighborhood of a point $\xi\in \mathbb{A}^n_k$. Let $E_i$ be a normal crossing divisor in $\mathbb{A}^n_k$. 
 
The composition ideal $K_i$ in $W_i$ of the product $J_i=M_i\cdot I_i$, with respect to a control $c_{i+1}$, is  
\begin{equation} \label{qus}
K_i= \left\{\begin{array}{ll} P_i\cdot I_{W_i}(E_i\cap W_i) & \text{ if }\ I_i\neq 1, \\ 1 & \text{
if }\ I_i=1. \end{array}\right.
\end{equation} 

The critical value for the following step of induction on the dimension is $c_i=ord_{\xi}(K_i)$.

The construction of the composition ideal $K_i$ ensures normal crossing with the exceptional divisor $E_i$. 

We say that an ideal $K$ is \emph{bold regular} if $K=<X^a>$, $K\in k[X]$, $a\in \mathbb{N}$.

Finally, construct the junior ideal $J_{i-1}$
\begin{equation} \label{jotas}
J_{i-1}= \left\{\begin{array}{ll} Coeff_V(K_i) & \text{ if }\ K_i \text{ is not bold regular or } 1 \\ 1 & \text{
otherwise }\ \end{array}\right.
\end{equation} where $V$ is a hypersurface of maximal contact in $W_i$ (see \cite{strong} page $830$) and $Coeff_V(K_i)$ is the coefficient ideal of $K_i$ in $V$ (see \cite{strong} page $829$). The junior ideal $J_{i-1}$ is an ideal in this suitable hypersurface $V$ of dimension $i-1$.

If $\frac{\theta_{n}}{c}\geq 1$ we are in the first case of equation 
(\ref{pes}), $\frac{\theta_{n-1}}{c_{n}}=\frac{\theta_{n-2}}{c_{n-1}}=\ldots
=\frac{\theta_j}{c_{j+1}}=1$ and $t_{j-1}=\ldots=t_1=\infty$ for $n-1\geq j\geq 1$, because $D_{n-1}=\ldots=D_1=\emptyset$ and $P_i=I_i$, and hence $J_{i-1}$
is always given by a unique monomial.

\begin{remark}
For an ideal $J=<X_1^{a_1}\cdot \ldots \cdot X_n^{a_n}>$ as in equation (\ref{jota}), if we assume $a_n\geq a_{n-1}\geq \ldots \geq a_1\geq c$, then at every stage $\frac{\theta_{n}}{c}\geq 1$
, so we are always in the previous situation. The singular locus of $(J,c)$ is always a union of hypersurfaces $\cup_{i=1}^r\{X_i=0\}$, $1\leq r \leq n$, and the center of the next monoidal transformation will be the intersection of some of these hypersurfaces. So we will call this case the {\bf minimal codimensional
case}. 
\end{remark}

\begin{remark} \label{remgrande}
If there exists some $a_{i_0}<c$, at a certain stage of
the resolution process it may occur $\frac{\theta_{n}}{c}< 1$. Then
we are in the second case of equation (\ref{pes}), the (exceptional) monomial part $M_n$
can appear in some $J_j$ for $n-1\geq j\geq 1$, and
$\frac{\theta_{j}}{c_{j+1}}$ can be much greater than $1$, what
increase the number of monoidal transformations. Now its singular locus is a union of intersections of hypersurfaces of the type $\cup_{l_j}(\{X_{l_1}=0\}\cap \ldots \cap \{X_{l_i}=0\})$. This is the {\bf higher codimensional case}.
\end{remark}

\section{Bound in the minimal codimensional case}

\begin{remark} From now on, we always look to the points where the function $t$, defined in (\ref{invt}), is maximal. So the following results concerning the behaviour of the function $t$ always affect the points where it reaches its maximal value. 
\end{remark}

\begin{proposition} \label{bajatn} Let $(W,(J,c),E)$ be a basic object where $J$ is a monomial ideal as in equation (\ref{jota}), with $a_i\geq c$ for all $1\leq i \leq n$. We can factor $J=J_n=M_n\cdot I_n$, and after $r-1$ monoidal transformations, $J_n^{(r-1)}=M_n^{(r-1)}\cdot I_n^{(r-1)}$. Let $\xi\in W^{(r-1)}$ be a point where  $ord_{\xi}(I_n^{(r-1)})=\theta_n$. After each monoidal transformation $\pi$, the resolution function in a neighbourhood of $\xi'$, where $\pi(\xi')=\xi$ and $ord_{\xi'}(I_n^{(r)})=\theta_n'<\theta_n$, is of the form $$\left(
\left[\frac{d-\sum_{j=1}^sa_{i_j}}{c},s\right],[1,0],\ldots,[1,0]\right)
\text{ for some }\  1\leq s\leq n-1. $$
\end{proposition}

\begin{proof}
After monoidal transformations, $$J_n^{(r)}=M_n^{(r)}\cdot I_n^{(r)}=<X_{i_1}^{b_1}\cdots
X_{i_s}^{b_s}>\cdot <X_{i_{s+1}}^{a_{i_{s+1}}}\cdots
X_{i_n}^{a_{i_n}}>$$  with
$d-\sum_{j=1}^sa_{i_j}=\sum_{j=s+1}^n a_{i_j}\geq c $ then,
$P_n^{(r)}=I_n^{(r)}$ and the (exceptional) monomial part does not appear in $J_l^{(r)}$
for all $n\geq l\geq 1$.

We have $\theta_n'\neq \theta_n$, then $E_n^{(r)}=Y'+|E|^{\curlyvee}$ and $m_n=s$, we count all the exceptional divisors of the previous steps and the new one. There are no exceptional divisors in lower dimension because $E_{n-1}^{(r)}=(Y'+|E|^{\curlyvee})-E_n^{(r)}=\emptyset$ and,
in a similar way, we obtain $E_l^{(r)}=\emptyset$ for all $n-1 \geq l\geq 1$.  

The normal crossing divisors $D_i^{(r)}=\emptyset$ for all $n-1\geq i \geq 1$, so the corresponding ideals $M_{n-1}^{(r)}=\ldots=M_1^{(r)}=1$. In particular, $M_{n-1}^{(r)}=1$, hence 
$$c_n'=ord_{\xi'}(K_n^{(r)})=ord_{\xi'}(Coeff(K_n^{(r)}))=ord_{\xi'}(J_{n-1}^{(r)})=ord_{\xi'}(I_{n-1}^{(r)})=\theta_{n-1}'$$ with $\xi'\in W^{(r)}$ such that $\pi(\xi')=\xi$, because $ord(Coeff(K))=ord(K)$ when $K$ is a monomial ideal, therefore $\frac{\theta_{n-1}'}{c_n'}=1$. By the same argument we obtain $\frac{\theta_{n-2}'}{c_{n-1}'}=\ldots=\frac{\theta_1'}{c_{2}'}=1$.
\end{proof}

\begin{remark} After each monoidal transformation, the exceptional divisors at each dimension are:
$$E_j'=\left\{ \begin{array}{ll} {\scriptstyle E_j^{\curlyvee}} &
\text{ \scriptsize if } {\scriptstyle
(t_n'(\xi'),\ldots,t_{j+1}'(\xi'))=(t_n(\xi),\ldots,t_{j+1}(\xi))} \text{ \scriptsize
and } {\scriptstyle \theta_j'=\theta_j} \\ {\scriptstyle
(Y'+(E_1\cup \ldots \cup E_n)^{\curlyvee})-(E_n'+\cdots
+E_{j+1}')} & \text{ \scriptsize in other case}
\end{array} \right.$$  \hspace*{0.5 cm}
${\scriptstyle ( E_n'=E_n^{\curlyvee}}$ {\scriptsize if}
${\scriptstyle \theta_n'=\theta_n}$ {\scriptsize or}
${\scriptstyle E_n'=Y'+(E_1\cup \ldots \cup E_n)^{\curlyvee}}$
{\scriptsize otherwise)} \vspace*{2 mm} \\ for $n> j\geq 1$, where $E_j^{\curlyvee}$
denotes the strict transform of $E_j$ by the monoidal transformation $\pi$, $Y'$ denotes
the new exceptional divisor, the point $\xi'\in W_i'$ satisfies $\pi(\xi')=\xi$, $\theta_j'=ord_{\xi'}(I_j')$ and $\theta_j=ord_{\xi}(I_j)$. We denote $|E|=E_1\cup \ldots \cup E_n$. 

Hence, after the first monoidal transformation,
since $\theta_n'<\theta_n$ we have $E_n'=Y'$ and
$E_{n-1}'=\cdots=E_1'=\emptyset$. After the second monoidal transformation, at the
chart where $\theta_n''=\theta_n'$ we obtain
$E_n''=(E_n')^{\curlyvee}=\emptyset$, $E_{n-1}''=Y''$, and
$E_{n-2}''=\cdots=E_1''=\emptyset$ and so on. We call this
phenomena propagation because every exceptional divisor appears in the resolution function $t$
firstly in dimension $n$, then in dimension $n-1$, $n-2$, and so on.
\end{remark}

\begin{definition} We will call {\bf propagation}, ${\bf p(i,j)}$, for $1\leq i\leq
j-1$, $1\leq j\leq n$, to the number of monoidal transformations needed to
eliminate $i$ exceptional divisors in dimension $j$, when we remain
constant $(t_n,t_{n-1},\ldots ,t_{j+1})$ and $\theta_j$, and there are
no exceptional divisors in lower dimensions $j-1,\ldots,1$. That is, passing from the stage  $$([\theta_n,m_n],\ldots,[\theta_{j+1},m_{j+1}],[\theta_j,i],[1,0],\ldots,
[1,0])$$ to the stage $$([\theta_n,m_n],\ldots,[\theta_{j+1},m_{j+1}],[\theta_j,0],[1,0],\ldots,
[1,0],\overbrace{\infty,\ldots,\infty}^{i}).$$
\end{definition}

\begin{lemma}{\bf Propagation Lemma} Let $(W,(J,c),E)$ be a basic object where $J$ is a monomial ideal as in equation (\ref{jota}) with $a_l\geq c$ for all $1\leq l \leq n$. Let $p(i,j)$ be the propagation of $i$ exceptional divisors in dimension $j$ in the resolution process of $(W,(J,c),E)$. \\ Then, for all $1\leq j\leq n$,
\begin{equation} \label{pro}
p\ (i,j)=\left\{ \begin{array}{ll} i+\sum_{k=1}^i p\ (k,j-1) &
\text{\ if \ } 0\leq i\leq j-1  \\ 0  & \text{\ if \ } i=j \\
\end{array} \right.\end{equation}
\end{lemma}

\begin{proof}
\begin{itemize}
\item If there are $i$ exceptional divisors in dimension $i$, $K_{i+1}$ is bold regular, $t_i=\infty$ then $p(i,i)=0$. We can not propagate these $i$ exceptional divisors at this stage of the resolution process. 

If there are $s$ exceptional divisors at this step of the resolution process, then there are $n-s$ variables in $I_n$. On the other hand, from dimension $n$ to dimension $i+1$ there are $s-i$ exceptional divisors.

When we construct $J_{n-1},\ldots ,J_{i+1}$, add to the corresponding composition ideal $K_j$ the variables in $I_{W_j}(E_j\cap W_j)$, so in these dimensions there are $(n-s)+(s-i)=n-i$ variables.

When we make induction on the dimension, at each step lose one variable, so in $n-i-1$ steps obtain that $K_{i+1}$, that corresponds to the $n-(n-i-1)=i+1$ position, is bold regular. And the variables appearing in these $i$ exceptional divisors do not appear in the center of the next monoidal transformation.   

\item By induction on the dimension:
\begin{itemize}
\item[-] If $j=1$, $p(1,1)=0$ by the previous argument.
\item[-] If $j=2$, $p(1,2)=1$ because when we propagate $1$ excepcional divisor from dimension
$2$ to dimension $1$, $K_2'$ is bold regular.

\begin{center} \vspace*{2 mm}
\begin{tabular}{c}
  $([\theta_n,m_n],\ldots,[\theta_2,1],[1,0])$  \\
  $\downarrow {\scriptstyle X_i}$ \\
  $([\theta_n,m_n],\ldots,[\theta_2,0],\infty)$ \\
\end{tabular} \vspace*{2 mm}
\end{center}
Then $p(1,2)=1=1+0=1+p(1,1)$.

\item[-]We assume that the result holds for $j\leq s-1$. For $j=s$:
$$\begin{array}{c} ([\theta_n,m_n],\ldots,[\theta_{s+1},m_{s+1}],[\theta_s,i],[1,0],\ldots,
[1,0])  \\ \downarrow  \\ 
([\theta_n,m_n],\ldots,[\theta_{s+1},m_{s+1}],[\theta_s,i-1],[1,1],[1,0],\ldots,[1,0])  \\ \vspace*{0.2cm}
\left. \hspace*{2 cm} \begin{array}{r} \downarrow  \\ \vdots
\\ \downarrow  \end{array} \right\} p(1,s-1) \\ \vspace*{0.2cm} ([\theta_n,m_n],\ldots,[\theta_{s+1},m_{s+1}],[\theta_s,i-1],[1,0],\ldots,
[1,0],\infty)  \\ \downarrow  \\ 
([\theta_n,m_n],\ldots,[\theta_{s+1},m_{s+1}],[\theta_s,i-2],[1,2],[1,0],\ldots,[1,0])
\\ \end{array}$$

We want $([\theta_n,m_n]\ldots [\theta_{s+1},m_{s+1}])$ and $\theta_s$ remain constant. So after the first monoidal transformation look to some suitable chart where $m_s=i$ drops. As $m_s$ drops then $m_{s-1}=i-(i-1)=1$ and propagate this exceptional divisor in dimension $s-1$, making $p(1,s-1)$ monoidal transformations. Otherwise, to keep $([\theta_n,m_n]\ldots [\theta_{s+1},m_{s+1}])$ and $\theta_s$ constant, the only possibility is to look to a suitable chart where $m_s$ drops from $i-1$ to $i-2$. But in this case this would provide the same resolution function that appears after the propagation. As we want to construct the largest possible sequence of monoidal transformations, we follow the propagation phenomenon as above.

After more monoidal transformations: 
$$\begin{array}{c} \left. \hspace*{2 cm} \begin{array}{r} \downarrow  \\ \vdots
\\ \downarrow  \end{array} \right\} p(2,s-1) \vspace*{0.2cm} \\ 
 \downarrow  \\ \vdots
\\ \downarrow  \\ ([\theta_n,m_n],\ldots,[\theta_{s+1},m_{s+1}],[\theta_s,1],[1,0],\ldots,
[1,0],\overbrace{\infty,\ldots,\infty}^{i-1}) \\  \downarrow  \\ 
([\theta_n,m_n],\ldots,[\theta_{s+1},m_{s+1}],[\theta_s,0],[1,i],[1,0],\ldots,[1,0])
\\ \end{array}$$ 
Then, for $1\leq i \leq s-1$,
$$p(i,s)=1+p(1,s-1)+1+p(2,s-1)+ \cdots +1+p(i,s-1)$$ with
$p(l,s-1)$, $1\leq l \leq i$, defined by the induction hypothesis. 
\end{itemize}
\end{itemize} 
\end{proof}

\begin{remark} Computation of examples in Singular with \emph{desing} package has been useful to state this behaviour of the exceptional divisors after monoidal transformations. The implementation of this package is based on the results appearing in \cite{paperlib}.  
\end{remark}

\begin{theorem} \label{inv} Let $(W,(J,c),\emptyset)$ be a basic object where $J$ is a monomial ideal as in equation (\ref{jota}) with $a_i\geq c$ for all $1\leq i \leq n$. Then, the resolution function corresponding to $(W,(J,c),\emptyset)$ drops after monoidal transformations in the following form:
\begin{center} 
\begin{tabular}{cl}
 \multicolumn{1}{c}{} &
 \multicolumn{1}{c}{} \\
$([\frac{d}{c},0],[1,0],\ldots,[1,0])$ & \\ $\hspace*{0.4 cm}
\downarrow {\scriptstyle X_i}$ & $1^{st}$  monoidal transformation \\
$([\frac{d-a_i}{c},1],[1,0],\ldots,[1,0])$ & \\ $ \hspace*{0.6 cm}
 \begin{array}{cc} \downarrow {\scriptstyle X_i} & \\ \vdots
\\ \downarrow {\scriptstyle X_i} & \\ \end{array} $ & $\left.
\hspace*{-7 mm} \begin{array}{c} \\ \\ \\ \end{array} \right\} p\
(1,n)$ monoidal transformations \\
$([\frac{d-a_i}{c},0],[1,0],\ldots,[1,0],\infty)$ & \\
$\downarrow$ & center defined only by variables in $I$ \\
$([\frac{d-a_i-a_j}{c},2],[1,0],\ldots,[1,0])$ & \\ 
$ \hspace*{0.3cm} \begin{array}{cc} \downarrow & \\ \vdots & \\ \downarrow &
\\ \end{array}$ & $\left. \hspace*{-7 mm} \begin{array}{c} \\ \\ \\
\end{array} \right\} p\ (2,n)$ monoidal transformations \\
$([\frac{d-a_i-a_j}{c},0],[1,0],\ldots,[1,0],\infty,\infty)$ & \\
$\downarrow$ & center defined only by variables in $I$ \\ $\vdots$ &
\hspace*{1 cm} $\vdots$ \\ $\downarrow$ & \\
$([\frac{a_l}{c},n-1],[1,0],\ldots,[1,0])$ & \\ $ \hspace*{0.3 cm}
\begin{array}{cc} \downarrow  & \\ \vdots & \\ \downarrow & \\
\end{array}$ & $\left. \hspace*{-7 mm}
\begin{array}{c} \\ \\ \\ \end{array} \right\} p\ (n-1,n)$ monoidal transformations \\
$([\frac{a_l}{c},0],\infty,\ldots,\infty)$ \hspace*{0.6 cm} & \\
\end{tabular}
\end{center}

At this stage, $a_l\geq c$ by hypothesis, so the center of the next monoidal transformation is $\{X_l=0\}$, and then we obtain an exceptional monomial.
\end{theorem}

\begin{proof}
It follows by the propagation lemma and the fact that each time that $\theta_n$ drops
$E_n'=Y'+|E|^{\nu}\neq\emptyset$, and $E_l'=(Y'+|E|^{\nu})-(E_n'+\cdots +E_{l+1}')=\emptyset$ for all
$n-1\geq l\geq 1$. 
\end{proof}

\begin{remark}
Following the propagation in the previous way provides the largest branch in the resolution tree, because in other case, for example after the first monoidal transformation
$$\begin{array}{c} ([\frac{d-a_i}{c},1],[1,0],\ldots,[1,0])  \\  
{\scriptstyle X_i}\swarrow \qquad \searrow {\scriptstyle X_j} \\ ([\frac{d-a_i}{c},0],[1,1],[1,0]\ldots,[1,0]) \qquad ([\frac{d-a_i-a_j}{c},2],[1,0],\ldots,[1,0])  \end{array}$$ looking to some chart $j$ with $j\neq i$
we obtain an invariant which will appear later in the resolution process, after the propagation $p(1,n)$.
\end{remark}

\begin{corollary}
Let $(W,(J,c),\emptyset)$ be a basic object where $J$ is a monomial ideal as in equation (\ref{jota}) with $a_i\geq c$ for all $1\leq i \leq n$. Therefore the number of monoidal transformations needed to transform $J$ into an exceptional
monomial is at most \begin{equation} \label{ec} 1 + p(1,n) +1 +
p(2,n)+ \ldots + 1 + p(n-1,n) + 1 =\  n + \sum_{j=1}^{n-1}p(j,n).
\end{equation}
\end{corollary}

\begin{remark} In this case we always have $\theta_n\geq c$, so
$Sing(J,c)\neq \emptyset$ at every stage of the resolution process. Therefore, in the resolution tree, the branch of
theorem \ref{inv} effectively appears, and it is the largest, hence (\ref{ec}) is exactly the number of monoidal transformations to obtain $J'=M'$.
\end{remark}

\begin{proposition} \label{cata} Let $(W,(J,c),\emptyset)$ be a basic object where $J$ is a monomial ideal as in equation (\ref{jota}) with $a_i\geq c$ for all $1\leq i \leq n$. Then the previous sum of propagations is a partial sum of Catalan numbers.  $$n +
\sum_{{\scriptscriptstyle j=1}}^{{\scriptscriptstyle n-1}}p(j,n)=
\sum_{{\scriptscriptstyle j=1}}^{{\scriptscriptstyle n}}C_j\
\text{  where }\ C_j=\left\{\frac{1}{j+1}\left(\begin{array}{c}
2j\\ j \end{array} \right)\right\}\ \text{ are Catalan numbers.}$$
\end{proposition}
\begin{proof} 
\begin{enumerate}
\item[(1)] Extend $p$ to arbitrary dimension: $$n + \sum_{j=1}^{n-1}p(j,n)=p(n,n+1).$$ Because of the form of the recurrence equation defining $p(i,j)$, and the fact that $p(n,n)=0$ by definition, it follows that \begin{equation} \label{rec} p(n,n+1)=n + \sum_{j=1}^{n}p(j,n)=n + \sum_{j=1}^{n-1}p(j,n). \end{equation} 
\item[(2)] Solve the recurrence equation defining $p(i,j)$:
\begin{enumerate}
\item We transform the recurrence equation (\ref{pro}), defining $p(i,j)$ for $0\leq i\leq j$ and $1\leq j\leq n$, to another recurrence equation defined for every $i,j\geq 0$:  

By sending the pair $(i,j)$ to the pair $(i,j-i)$ we extend the recurrence to $i,j\geq 0$, that is, we consider $$\tilde{p}(i,j)=p(i,i+j)$$ then $p(i,j)=\tilde{p}(i,j-i)$. As $p(i,j)$ is defined for $0\leq i \leq j$ then $\tilde{p}(i,j)$ is defined for $0\leq i \leq i+j$ for every $i,j\geq 0$. 

\item Note that $$\tilde{p}(i,j)-\tilde{p}(i-1,j+1)=p(i,i+j)-p(i-1,i+j)$$ $$=i+\sum_{k=1}^i p(k,i+j-1)-(i-1)- \sum_{k=1}^{i-1}p(k,i+j-1)=p(i,i+j-1)+1=\tilde{p}(i,j-1)+1.$$ Therefore, we have the following recurrence equation involving $\tilde{p}(i,j)$ 
\begin{equation} \label{chis} \left\{\begin{array}{ll}\tilde{p}(i,j)=1+\tilde{p}(i-1,j+1)+\tilde{p}(i,j-1) & \text{ for } i,j \geq 1 \\ \tilde{p}(0,j)=\tilde{p}(i,0)=0 & \end{array}\right. \end{equation}
Take $r(i,j)=p(i,i+j)+1=\tilde{p}(i,j)+1$ and replace $\tilde{p}(i,j)$ with $r(i,j)$ in the equation (\ref{chis}). It follows the auxiliary recurrence equation:
\begin{equation} \label{erres} \left\{\begin{array}{ll} r(i,j)=r(i-1,j+1)+r(i,j-1)& \text{ for } i,j \geq 1\\ r(0,j)=\tilde{p}(0,j)+1=1,\   r(i,0)=\tilde{p}(i,0)+1=1 & \end{array}\right. \end{equation}

\item Resolving the auxiliary recurrence equation (\ref{erres}) by generating functions: 

Define $r_{i,j}:=r(i,j)$ and the generating functions $$R(x,y)=\sum_{i,j\geq 0}r_{i,j}x^iy^j \in \mathbb{C}[[x,y]], \   R_s(x,y)=\sum_{i,j\geq 1}r_{i,j}x^{i-1}y^{j-1} \in \mathbb{C}[[x,y]].$$ Note that $R(x,y)$ is, by definition, the generating function of the sequence $r(i,j)$.

By the recurrence equation (\ref{erres}) involving $r(i,j)$, it follows 
$$\hspace*{-3.5cm} R_s(x,y)=\sum_{i,j\geq 1}r_{i-1,j+1}x^{i-1}y^{j-1}+\sum_{i,j\geq 1}r_{i,j-1}x^{i-1}y^{j-1}$$
$$ \hspace*{1cm} \begin{array}{l}
=\sum\limits_{i\geq 0,j\geq 1}r_{i,j+1}x^{i}y^{j-1}+\frac{1}{x}\sum\limits_{i\geq 1,j\geq 0}r_{i,j}x^{i}y^{j} \vspace*{0.15cm} \\ 
=\frac{1}{y^2}\sum\limits_{i\geq 0,j\geq 1}r_{i,j+1}x^{i}y^{j+1}+\frac{1}{x}\left[\sum\limits_{i\geq 1}r_{i,0}x^{i}+ \sum\limits_{i\geq 1,j\geq 1}r_{i,j}x^{i}y^{j}\right]  \vspace*{0.15cm} \\
=\frac{1}{y^2}\sum\limits_{i\geq 0,j\geq 2}r_{i,j}x^{i}y^{j}+\frac{1}{x}\left[\sum\limits_{i\geq 1}x^{i}+ \sum\limits_{i\geq 1,j\geq 1}r_{i,j}x^{i}y^{j}\right]  \vspace*{0.15cm} \\
=\frac{1}{y^2}\left[\sum\limits_{i\geq 0,j\geq 1}r_{i,j}x^{i}y^{j}-\sum\limits_{i\geq 0}r_{i,1}x^{i}y\right] +\frac{1}{x}\left[\frac{1}{1-x}-1+ xyR_s(x,y)\right]  \vspace*{0.15cm} \\
=\frac{1}{y^2}\left[\sum\limits_{j\geq 1}r_{0,j}y^{j}+xyR_s(x,y)-y\sum\limits_{i\geq 0}r_{i,1}x^{i}\right] +\frac{1}{x}\left[\frac{x}{1-x}+ xyR_s(x,y)\right] \vspace*{0.15cm} \\
=\frac{1}{y^2}\left[\frac{y}{1-y}+xyR_s(x,y)-y\sum\limits_{i\geq 0}r_{i,1}x^{i}\right] +\frac{1}{1-x}+ yR_s(x,y) \vspace*{0.15cm} \\
=\frac{1}{y(1-y)}+\frac{x}{y}R_s(x,y)-\frac{1}{y}\sum\limits_{i\geq 0}r_{i,1}x^{i}+\frac{1}{1-x}+ yR_s(x,y). 
\end{array}$$  

Then $$\left(1-y-\frac{x}{y}\right)R_s(x,y)=\frac{1}{y(1-y)}+\frac{1}{1-x}-\frac{1}{y}\sum_{i\geq 0}r_{i,1}x^{i}$$ multiplying the equality by \emph{y} we have $$(y-y^2-x)R_s(x,y)=\frac{1}{1-y}+\frac{y}{1-x}-\sum_{i\geq 0}r_{i,1}x^{i}$$
$$=\frac{1}{1-y}+\frac{y}{1-x}-r_{0,1}-\sum_{i\geq 1}r_{i,1}x^{i}=\frac{y}{1-y}+\frac{y}{1-x}-\sum_{i\geq 1}r_{i,1}x^{i}.$$ Therefore $$(y-y^2-x)R_s(x,y)=\frac{y}{1-y}+\frac{y}{1-x}-\sum_{i\geq 1}r_{i,1}x^{i}$$ which defines an equation of the form $$Q(x,y)R_s(x,y)=K(x,y)-U(x).$$ Now apply the \emph{kernel method} used in \cite{pet}, algebraic case $4.3$: 

If $Q(x,y)=0$ then $y=\frac{1\pm \sqrt{1-4x}}{2}$. We take the solution passing through the origin, $y=\frac{1- \sqrt{1-4x}}{2}$ and $y=xC(x)$ where $C(x)$ is the generating function of Catalan numbers. 

On the other hand, $Q(x,y)=0$ gives $K(x,xC(x))=U(x)$, $$K(x,y)=\frac{y}{1-y}+\frac{y}{1-x}=\frac{-y^2+y-x+1}{(1-x)(1-y)}-1$$ so $K(x,xC(x))=\frac{1}{(1-x)(1-xC(x))}-1$ and using $\frac{1}{1-xC(x)}=C(x)$ we have $$U(x)=\frac{C(x)}{1-x}-1.$$ Making some calculations and using that $R(x,y)$ satisfies $$R(x,y)=xyR_s(x,y)+\sum_{j\geq 0}r_{0,j}y^j+\sum_{i\geq 0}r_{i,0}x^i-r_{0,0}$$ we obtain the generating function of $r(i,j)$ $$R(x,y)=\frac{xyC(x)+x-y}{(y^2-y+x)(1-x)}.$$
\end{enumerate}

\item[(3)] Compute the generating function of the sequence $p(n,n+1)$:

The coefficient of \emph{y} in $R(x,y)$ is just $\sum_{i\geq 0}r_{i,1}x^{i}$ then $$\sum_{i\geq 0}r_{i,1}x^{i}=\frac{\partial R(x,y)}{\partial y}\Big|_{y=0}= \frac{C(x)}{1-x}$$ is the generating function of the elements in the first column. 

If $C(x)$ is the generating function of $C_n$ then the convolution product $C(x)\cdot \frac{1}{1-x}$ is the genera\-ting function of $\sum_{k=0}^nC_k=S_n$ therefore $$r_{n,1}=\sum_{k=0}^nC_k.$$ As $r(n,1)=p(n,n+1)+1$ then $p(n,n+1)=r(n,1)-1=\sum_{k=0}^nC_k-1$, as $C_0=1$ we have $$p(n,n+1)=\sum_{k=1}^nC_k$$ where $C_k$ are Catalan numbers. 
\end{enumerate}  
\end{proof}

See \cite{stan} for more details about Catalan numbers and the web
page \cite{slo} for further details about their partial sums.

\begin{theorem} Let $(W,(J,c),\emptyset)$ be a basic object where $J$ is a monomial ideal as in equation (\ref{jota}) with $a_i\geq c$ for all $1\leq i \leq n$. Then the number of monoidal transformations required to resolve $(W,(J,c),\emptyset)$ is at most $$\sum_{j=1}^{n}C_j+(2^{\sum_{j=1}^{n}C_j}-1)(d-c)-c+1$$ where $C_j$ are Catalan numbers.
\end{theorem}

\begin{proof} It follows by theorem \ref{mon}, lemma \ref{grad} and  proposition \ref{cata}.
\end{proof}

\begin{example}
The following table shows some values of the bound for any monomial ideal $J$ as in equation (\ref{jota}) with $a_i\geq c$ for all $1\leq i \leq n$.

\vspace*{0.6cm}

\begin{table}[ht] \vspace*{-0.6cm}
\caption{Values of the bound}
\begin{center} 
\begin{tabular}{||c|c|c||} \hline 
 {$n$}  & $\sum_{j=1}^{n}C_j$ & global bound \\ \hline
 $1$ & $1$ & $1+(d-c)-c+1$ \\ \hline
 $2$ & $3$ & $3+7(d-c)-c+1$ \\ \hline
 $3$ & $8$ & $8+255(d-c)-c+1$ \\ \hline
 $4$ & $22$ & $22+4194303(d-c)-c+1$ \\ \hline
\end{tabular}
\end{center}
\end{table}
\end{example}

\vspace*{-0.5cm}

\begin{remark} Note that, as a consecuence of proposition \ref{cata}, the number of monoidal transformations needed to transform $J$ into an exceptional monomial only depends on $n$, the dimension of the ambient space.
\end{remark}

\begin{corollary} Let $J=<Z^c-X_1^{a_1}\cdot \ldots \cdot X_n^{a_n}>\subset k[X_1,\ldots,X_n,Z]$ be a toric ideal with $a_i\geq c$ for all $1\leq i \leq n$. Then the number of monoidal transformations needed to resolve $(\mathbb{A}^{n+1}_k,(J,c),\emptyset)$ is at most $$\sum_{j=1}^{n}C_j+(2^{\sum_{j=1}^{n}C_j}-1)(d-c)-c+1 $$ where $C_j$ are Catalan numbers and $d=\sum_{i=1}^n a_i$.
\end{corollary}

\section{Higher codimensional case}

In the minimal codimensional case, the way in which the invariant drops essentially
depends on the number of accumulated exceptional divisors. Because the first components of the invariant, $\theta_{n},\ldots,\theta_{1}$, defined in equation (\ref{invt}), only depend on the order of the ideals $I_n,\ldots,I_1$. Recall that, for each $J_i$, we use the ideal $M_i$ (see remark \ref{constructM}), to define the ideal $I_i$.  

But in this case the first components of the invariant play an important
role. They can also depend on the order of the (exceptional) monomial
part $M_n$, see remark \ref{remgrande}. So they may increase suddenly when some $\theta_{j}$ is given by the order of the ideal $M_n$. We will call this situation the {\it higher codimensional case in dimension $j$}.

Note that after some monoidal transformations, we can obtain a new higher codimensional case in another dimension. 

So, we must compute the number of monoidal transformations while $\theta_n\geq c$, with a
suitable sum of propagations. Then, estimate the number of monoidal transformations needed to get the higher codimensional case in dimension $1$, and use the known estimation for the
order of $M_n$ to give an upper bound for the number of monoidal transformations needed to get the
following higher codimensional case inside this one (if it is possible). Afterward, estimate the number of monoidal transformations needed to get the higher codimensional case in dimension $2$,
and so on.

Hence, it has not been possible to obtain a bound for this case in
the same way as above, due to the complications of the
combinatorial problem, that perform that we can not know what
branch is the largest in the resolution tree (to obtain an
exceptional monomial).

Furthermore, if we could find such bound, the large number of potential cases we expect, suggests that this bound would be very huge, even to estimate only the number of monoidal transformations needed to obtain an exceptional
monomial.

\begin{example} If we consider the basic object $(W,(J,c),\emptyset)=(\mathbb{A}^3_k,(X_1^5X_2^4X_3,4),\emptyset)$, there exists a branch of height $15$ to obtain $J'=M'$ or
$Sing(J',c)=\emptyset$. So, in dimension $3$, we need a bound
greater than or equal to $15$ for a higher codimensional case, in front of the $8$ monoidal transformations needed for a minimal codimensional case.
\end{example}

\begin{remark}
In any case, both theorem \ref{mon} and lemma \ref{grad} are valid also in the higher codimensional case. So the open problem is to find a bound $C$ to obtain an
exceptional monomial, to construct a global bound of the form
$$C+(2^C-1)(d-c)-c+1 .$$
\end{remark}

\begin{remark}
For $n=2$ the higher codimensional case appears only in dimension $1$ and making
some calculations we obtain $C=3$, that gives the same bound as in
the minimal codimensional case. This bound can be improved by studying the
different branches.
\end{remark}

\section*{Acknowledgments}

I am very grateful to Professor O. Villamayor for his suggestions to improve the writing of this paper, and Professor S. Encinas for his continued advice and help. I also would like to thank the anonymous referee for useful comments to improve the presentation of the paper.

\noindent Universidad de Valladolid \\ Departamento de Matem\'atica Aplicada \\ E.T.S. Arquitectura, Avda. de Salamanca s/n. \\ $47014$ Valladolid Spain


\begin{thebibliography}{0}

\bibitem[1]{bm} E. Bierstone and P. Milman: {\it Canonical desingularization in characteristic zero by blo\-wing up the maximum strata of a local invariant}. Invent. math. ${\bf 128}$ $(1997)$, $207-302$.

\bibitem[2]{pet} M. Bousquet-M$\acute{\rm e}$lou and M.
Petkovsek: {\it Linear recurrences with constant coefficients: the
multivariate case}. Discrete Mathematics, ${\bf 225}$ ($2000$),
$1-3$, $51-75$.

\bibitem[3]{paperlib} G. Bodn\'ar and J. Schicho: {\it Automated resolution of singularities for hypersurfaces}. Journal of Symbolic Computation, ${\bf 30},\ 4\ (2000),\ 401-428$.

\bibitem[4]{lib} G. Bodn\'ar and J. Schicho: {\it ``desing''}- A computer program for resolution of singularities. {\it http://www.risc.uni-linz.ac.at/projects/basic/adjoints/blowup/}.

\bibitem[5]{strong} S. Encinas and H. Hauser: {\it Strong
resolution of singularities in characteristic zero}. Commentarii
Mathematici Helvetici, ${\bf 77},\ 4\ (2002),\ 821-845$.

\bibitem[6]{cour} S. Encinas and O. Villamayor: {\it
A Course on Constructive Desingularization and Equiva\-riance}. In
{\it Resolution of Singularities, A research textbook in tribute
to Oscar Zariski} (Basel, $2000$), H. Hauser, J. Lipman, F. Oort,
and A. Quir\'os, Eds. Progress in Math. ${\bf 181}$,
Birkh$\ddot{\rm a}$user, $147-227$.

\bibitem[7]{EncinasVillamayor2003} S. Encinas and O. Villamayor. {\it A new proof of desingularization
over fields of characteristic zero}. Rev. Mat. Iberoamericana,
${\bf 19(2)}$: $339-353$, $2003$.

\bibitem[8]{sing} G.-M. Greuel, G. Pfister, and H. Sch\"onemann.
{\sc Singular} 2.0.5. A Computer Algebra System for Polynomial
Computations. Centre for Computer Algebra, University of
Kaiserslautern (2003). {\it http://www.singular.uni-kl.de}.

\bibitem[9]{existence} H. Hironaka: {\it Resolution of singularities of an algebraic variety over a field of characteristic zero}. Annals of Mathematics, ${\bf 79}$ $(1964)$, $109-362$. 

\bibitem[10]{ko} J. Koll$\acute{\rm a}$r: {\it Lectures on Resolution of Singularities}. Annals of Mathematics Studies, ${\bf 166}$ (Princeton and Oxford, $2007$), Princeton University Press.

\bibitem[11]{slo} N. J. A. Sloane: Sequences number A$014138$ and A$000108$
from The On-Line Encyclopedia of Integer Sequences, {\it
http://www.research.att.com/$\sim$njas/sequences/index.html}.

\bibitem[12]{stan} R. P. Stanley: {\it Enumerative Combinatorics}.
Vol $2$. Cambridge Studies in Avanced Mathematics, ${\bf 62}$
(Cambridge, $1999$), Cambridge University Press.

\bibitem[13]{villa} O. Villamayor: {\it Constructiveness of Hironaka's resolution}. 
Annales Scientifiques \'Ecole Normale Sup\'erieure, $4^{e}$ serie, ${\bf 22}$ $(1989)$, $1-32$.

\bibitem[14]{Wo} J. Wodarczyk: {\it Simple Hironaka resolution in characteristic zero}. 
Journal of the Ameri\-cal Mathe\-matical Society, ${\bf 18}$ ($2005$), $779-822$. 


\end{thebibliography}
\end{document}